\documentclass[12pt]{article}
\usepackage{amsmath,amssymb,amsfonts,mathtools}
\usepackage{longtable}

\allowdisplaybreaks

\def\Z{{\mathbb Z}}

\def\R{{\mathbb R}}
\def\C{{\mathbb C}}

\usepackage{color}

\title{A Remark on Fourier Transform}
\author{Guangwu Xu\thanks{SCST, Shandong University,
Qingdao, Shandong, China;
e-mail: {\tt gxu4sdq@sdu.edu.cn}.}
}
\date{}

\begin{document}
\maketitle

\begin{abstract} In  this note, we describe an interpretation of the (continuous) Fourier
transform from the perspective of the Chinese Remainder Theorem. Some related issues, including a new derivation 
of Poisson summation formula, are discussed.
\end{abstract}


\section{Introduction}
The usual Chinese Remainder Theorem (CRT) concerns  reconstructing an integer given its remainders with respect to
a set of coprime moduli. For a set of pairwise coprime positive integers $\{m_0, m_1, \dots , m_{\mu -1}\}$,  we let $ \Gamma=\prod_{j=0}^{\mu -1} m_j$.
Then there is an isomorphism
\begin{eqnarray}
\Z/(\Gamma) &\rightarrow& \prod_{j=0}^{\mu -1}\Z/(m_j)\label{iso}\\
n &\mapsto& \bigg( n \pmod{m_0}, n \pmod{m_1},\cdots, n \pmod{m_{\mu -1}}\bigg).\notag
\end{eqnarray}
The crucial part is to verify the inverse (of the isomorphism), which is given by the CRT:
 finding (the unique) solution $n< \Gamma$ of the system of congruence equations
\begin{equation}\label{eq:crt}
\left\{ \begin{array}{l} x\equiv r_0 \pmod {m_0} \\
                               x\equiv r_1 \pmod {m_1} \\
                                \cdots \\
                               x\equiv r_{\mu -1} \pmod {m_{\mu -1}}
\end{array} \right.
\end{equation}
CRT provides a formula  (reported by Jiushao Qin (aka Chin Chiu-shao)  in 1247 \cite{Qin})
for the solution:
\begin{equation}\label{eq:1.1}
n = \sum_{j=0}^{\mu -1} r_j u_j \frac{\Gamma}{m_j} \pmod {\Gamma}
\end{equation}
where $u_j = \big(\frac{\Gamma}{m_j}\big)^{-1} \pmod {m_j}$. It is remarked that in \cite{Qin}, Qin explained a beautiful algorithm
for computing (the positive value) $u_j$ (under the name of `Dayan deriving one'). As discussed in \cite{Xu,XuLi}, Qin's algorithm
is somehow more concise and efficient than the modern Extended Euclidean algorithm.

We would like to point out that Qin also noted the following relation
\begin{equation}\label{eq:use}
 \sum_{j=0}^{\mu -1}  u_j \frac{\Gamma}{m_j} = 1 + \ell \Gamma,
\end{equation}
and gave it the names `positive use' (for $\ell=1$) and `universal use' (for $\ell>1$).
Using this relation to get (\ref{eq:1.1}) is more convenient as illustrated in \cite{dlx}.

Modeling a computational task to the Chinese remainder representation is of great significance. Each component of
the representation is restricted in a ring (e.g., $\Z/(m_j)$) of smaller size and
arithmetic operations on those smaller rings
are closed and independent, so computation task of large scale can be
decomposed into that of smaller scales and
performed in parallel, see  \cite{bch,dlx}.
Other applications of CRT in fast computation can be found
in public key cryptography.

It is well-known that the Chinese remainder theorem
 has a general formulation to decompose a certain ring to be a product of `smaller' quotient rings. The correctness of the
 decomposition  in essence is (\ref{eq:1.1}).
As an example, replacing  $\Z/(\Gamma)$ with $\R[x]/(x^n-1)$ and  $\Z/(m_j)$ with $\R[x]/(x-e^{\frac{-2j\pi i}n})$, we get the (finite) discrete Fourier transform: for a polynomial $f\in \R[x]$ with degree less than $n$,
\[
 f\mapsto \bigg(f(1), f(e^{\frac{-2\pi i}n}), f(e^{\frac{-4\pi i}n}), \cdots, f(e^{\frac{-2(n-1)\pi i}n})\bigg).
\]
The inverse of the transform is exactly in the same format like (\ref{eq:1.1}) and it is  now commonly called the Lagrange interpolation formula.
For the computational significance of the  discrete Fourier transform,  besides being parallelizable,
the symmetry of the roots of $x^n-1=0$ enables the use of divide and conquer strategy, so we
have the famous Cooley-Tukey algorithm --- the fast Fourier transform (FFT).

In this note, we would like to remark that the continuous version of Fourier transform can also be interpreted as some
(approximation) form
of the Chinese remainder representation. Such algebraic treatment enables different and maybe elementary approach to
some related issues.

\section{Continuous Fourier Transformand Chinese Remainder Theorem}
Now we consider a discretization of the continuous Fourier transform under the framework of Chinese remainder theorem.
We shall assume functions involved are all good enough (e.g., in Schwartz space \cite{SS}) so the convergence will not be
an issue.

Let $f:\R \rightarrow \C$ be a  rapidly decreasing smooth function, its Fourier transform is given by
\begin{equation}\label{eq:2.1}
\hat f (y) = \int_{\R} f(x) e^{-2\pi i xy} dx
\end{equation}
The inverse of the transform is proved to be
\begin{equation}\label{eq:2.2}
f (x) = \int_{\R}\hat f  (y) e^{2\pi i xy} dy.
\end{equation}

We will use carefully formulated finite sums to approximate these two absolutely convergent integrals.
Let $N,M$ be large even integers (depend on $f$).

For the spacial domain $\R$, we choose interval $[-\frac{N}2,  \frac{N}2]$ and partition it with each unit interval being divided into $M$ equal
subintervals. Take the following partition points $x_j=-\frac{N}2+\frac{j}{M}$ for $j=0,1,\cdots, MN-1$.

For the frequency domain $\R$, we choose interval $[-\frac{M}2,  \frac{M}2]$ and partition it with each unit interval being divided into $N$ equal
subintervals.  Take the following partition points $y_k=-\frac{M}2+\frac{k}{N}$ for $k=0,1,\cdots, MN$-1.

To approximately evaluate the transform at partition points $y_k$, we see that
\begin{eqnarray} \label{eq:2.3}
\hat f (y_k) &\approx& \int_{-\frac{N}2}^{\frac{N}2} f(x) e^{-2\pi i xy_k} dx\approx \sum_{j=0}^{MN-1} f(x_j) e^{-2\pi i x_jy_k}\frac{1}M \notag\\
& = & \frac{1}M\sum_{j=0}^{MN-1} f(x_j) e^{-2\pi i \frac{(Mx_j)(Ny_k)}{MN}} .
\end{eqnarray}

Similarly, we can get
\begin{equation}\label{eq:2.4}
 f(x_j) \approx \frac{1}N  \sum_{k=0}^{MN-1} \hat f(y_k) e^{2\pi i \frac{(Ny_k)(Mx_j)}{MN}}
\end{equation}

We should note that $Ny_k$ (resp. $Mx_j$) runs over all integers from $-\frac{MN}2$ to $\frac{MN}2 -1$ for $k=0, 1, \cdots, MN-1$ (for $j=0, 1, \cdots, MN-1$), with the same order.
Write $w=\frac{NM}2$ and and $\omega = e^{-2\pi i \frac{1}{MN}}$. Denote
\scriptsize{
\[
P(X)= \frac{1}M\big( f(x_w)+f(x_{w+1})X+\cdots+f(x_{2w-1})X^{w-1}+
f(x_0)X^w+f(x_1)X^{w+1}+\cdots +f(x_{w-1})X^{2w-1} \big).
\]
}\normalsize{}
Then
(\ref{eq:2.3}) says that
\[
P(X)\pmod{X-\omega^{Ny_k}}\approx \hat f (y_k) = \hat f (y_{Ny_k+\frac{MN}2}),
\]
and hence
$\hat f$ can be approximated by the  Chinese remainder representation of $P(X)$ ( treated as an element of
$\C[X]/(X^{MN}-1)$  ).

With this interpretation, we can illustrate a way of recovering $f$ given the values of $\hat f (y_k), \ k = 0, 1, \cdots, NM-1$
by using (\ref{eq:1.1}). Note that in this case $\Gamma=X^{MN}-1$ and $m_j=X-\omega^j$.  We also note that
$\left(\frac{\Gamma}{m_j} (\omega^j)\right)^{-1}$ is $\left(\frac{\Gamma}{m_j}\right)^{-1} \pmod {m_j}$, i.e., $u_j = \frac{1}{\prod_{\ell=0, \ell \neq j}^{NM-1} (\omega^j-\omega^{\ell})}=\frac{\omega^j}{MN}$. We have heuristically that \footnote{It is interesting to note that we can omit
the $\pmod$ operation here since the degree
of both sides are less than NM. We also note that
$\sum_{j=0}^{NM-1}  u_j \frac{\Gamma}{m_j} = \sum_{j=0}^{NM-1} \frac{ \frac{\Gamma}{m_j}}{\frac{\Gamma}{m_j}(\omega_j)} = 1$. }
\[
P(X) \approx \sum_{k=0}^{NM-1} \hat f (y_{(k+\frac{MN}2)\pmod {NM}}) u_k \frac{\Gamma}{m_k}.
\]
We note that $\frac{1}M f(x_j)$ is the coefficient of $X^q$ in $P(X)$ with $q = (j + \frac{MN}2 ) \pmod { MN}$, in order to get its expression,
we take $q$th derivatives of both sides and evaluate them at $X = 0$:
\[
 f(x_j) \approx \frac{1}N  \sum_{k=0}^{MN-1} \hat f(y_k) e^{2\pi i \frac{(Ny_k)(Mx_j)}{MN}}.
\]

As mentioned above, this derivation of the approximation formula for inverse Fourier
transform we just illustrated is rather heuristic. A more rigorous treatment inside this
framework can be done by using  the Dirichlet kernel which is defined as
\[
D(x) =\left\{ \begin{array}{ll}  \frac{\sin(\pi Mx+\frac{\pi x}N)}{\sin (\frac{\pi x}N)} &  \mbox{ if } \frac{x}N\notin \Z\\
                              MN +1                    &  \mbox{ if } \frac{x}N\in \Z\\
\end{array}
\right.
\]
We shall simply explain this for the case of $x_j=0$, i.e., $j=\frac{MN}2$.
\begin{eqnarray*}
\frac{1}N  \sum_{k=0}^{MN-1} \hat f(y_k)&=& \frac{1}N  \sum_{k=0}^{MN-1}  \int_{\R} f(x) e^{-2\pi i xy_k} dx
=\frac{1}N \int_{\R} f(x) \sum_{k=0}^{MN-1} e^{-2\pi i xy_k} dx\\
&=&\frac{1}N \int_{\R} f(x) \sum_{k=0}^{MN} e^{-2\pi i \frac{x}N \left(-\frac{MN}2+k\right)} dx-
\frac{1}N \int_{\R} f(x)  e^{-2\pi i Mx} dx\\
&=&\frac{1}N \int_{\R} f(x)   D(x) dx+\frac{O(1)}N
\approx f(0).
\end{eqnarray*}

We also want to remark that parallelization and FFT are thus possible for numerical evaluations of the transform.

\subsection{Poisson summation formula}
The classical $1$-dimensional  Poisson summation formula states that for
a rapidly decreasing smooth function $f$, one has
\begin{equation}\label{eq:2.5}
\sum_{n\in \Z} f(n) = \sum_{n\in \Z} \hat{f}(n).
\end{equation}

 It is noted that from our previous formula (\ref{eq:2.1}),
\begin{eqnarray*}
\hat f(-\frac{M}2)&+&\hat f(-\frac{M}2+1)+\cdots+\hat f(\frac{M}2-1)\\
&=&\hat f(y_0)+\hat f(y_N)+\hat f(y_{2N})+\cdots+\hat f(y_{(M-1)N})\\
&\approx & \frac{1}M \sum_{\ell=0}^{M-1}\sum_{j=0}^{MN-1}f(x_j)e^{-2\pi i \frac{(\frac{-MN}2+\ell N)Mx_j}{MN}}
       = \frac{1}M \sum_{j=0}^{MN-1}f(x_j) \sum_{\ell=0}^{M-1} \left(e^{-2\pi i \frac{ j}{M}}\right)^{\ell}\\
 &=&f(-\frac{N}2)+f(-\frac{N}2+1)+\dots+f(\frac{N}2-1)
 \end{eqnarray*}

Taking limits, we get the
Poisson summation formula (\ref{eq:2.5}).
This argument is essentially of the form of the finite group version of the Poisson summation formula.
This can be done because our partitions for approximating
integrals have a natural subgroup
structure.

\section*{Acknowledgement} The author thanks Dr. Shangbin Cui for useful discussions.


\begin{thebibliography}{9}

\bibitem{bch} P. Beame, S. Cook and J. Hoover, Log depth circuits for division and
related problems {\sl SIAM J. Comput.}, {\bf 15}:994--1003 (1986).
\bibitem{dlx} G. Davida, B. Litow , and G. Xu, Fast arithmetics using Chinese Remaindering, {\sl Information Processing
Letters}, 109(2009), 660-662.


\bibitem{Libbrect} U. Libbrect,
Chinese Mathematics in the Thirteenth Century, {\sl Dover Publications}, 2005.

\bibitem{Qin} J. Qin, Mathematical Treatise in Nine Sections, 1247.

\bibitem{SS} E. Stein and R. Shakarchi, Fourier Analysis--An Introduction, {\sl Princeton University Pres}, 2003.

\bibitem{Xu} G. Xu, On solving a generalized Chinese Remainder Theorem in the presence of remainder errors, {\sl Springer Proceedings in Math. \& Stat. Series ---Geometry, Algebra, Number Theory, and Their Information Technology Applications}, 2018 (to appear).
\bibitem{XuLi} G. Xu and B. Li, On the algorithmic significance and analysis of the method of DaYan Deriving One ({\sl in Chinese}),
{\sl https://arxiv.org/abs/1610.01175}.
\end{thebibliography}
\end{document}